\documentclass[11pt]{amsart}
\usepackage[latin1]{inputenc}

\usepackage{amsthm}
\usepackage{mathrsfs}
\usepackage{amsmath}
\usepackage{amssymb}

\theoremstyle{plain}

\newcommand{\e}{\varepsilon}

\newcommand{\oo}[1]{\mathscr O^{(#1)}}
\def\xx#1 {\newtheorem{#1}[thm]{#1}}

\xx Theorem
\xx Corollary
\xx Lemma

\theoremstyle{definition}

\xx Definition
\xx Example
\xx Examples
\xx Notation
\xx Proposition
\xx Remark
\xx Convention
\xx Conclusion
\xx Fact

\newcommand{\NN}{\mathbb N}
\newcommand{\RR}{\mathbb R}
\newcommand{\on}{{\upharpoonright}}

\DeclareMathOperator{\Pol}{Pol}

\author{Martin Goldstern}
\address{Algebra/DMG, TU Wien, Wiedner Hauptstr 8-10/104, 1040 Wien, Austria}
\urladdr{http://www.tuwien.ac.at/goldstern/}
\email{goldstern@tuwien.ac.at}
\thanks{The first author is supported by the Austrian Science Foundation FWF, grant  P 22994-N18.}

\author{G\'abor S\'agi}
\thanks{ The second author is supported by Hungarian
National Foundation for Scientific Research grant K68262 and by the J\'anos Bolyai Research Scholarship of the Hungarian Academy of
Sciences.}
\address{BUTE Department of Algebra and Alfr\'ed R\'enyi Institute of Mathematics, Hungarian Academy of Sciences,
Re\'altanoda u. 13-15, 1053 Budapest, Hungary}
\email{sagi@renyi.hu}

\author{Saharon Shelah}
\address{Einstein Institute of Mathematics\\
Edmond J. Safra Campus, Givat Ram\\
The Hebrew University of Jerusalem\\
Jerusalem, 91904, Israel\\
and
Department of Mathematics\\
Rutgers University\\
New Brunswick, NJ 08854, USA}%{4}
\thanks{The third author is supported by the German-Israeli Foundation for Scientific Research \& Development
   Grant No. 963-98.6/2007. Publication 989.}
\email{shelah@math.huji.ac.il}
\urladdr{http://shelah.logic.at/}

\title{Clones above the unary clone}
\begin{document}
\thispagestyle{plain}

\begin{abstract}
Let $\mathfrak c:= 2^{\aleph_0}$. We give a family of pairwise
incomparable clones on $\NN$ with $2^{\mathfrak c}$ members, all with the
same unary fragment, namely the set of all unary operations.

We also give, for each $n$, a family of $2^{\mathfrak c}$ clones all
with the same $n$-ary fragment, and all containing the set of all
unary operations.
\end{abstract}

\maketitle

\section{Introduction}

In this paper, $X$ will always be a countably infinite set. For a
fixed base set $X$, an operation on $X$ is a function $f:X^n\to
X$ for some positive natural number $n$. A clone on $X$ is a set
of operations that contains all projection functions and is
closed under composition. The set of all clones on $X$ ordered by
inclusion forms a complete lattice.  (The survey paper \cite{survey}
gives some background about clones, and in particular collects
many recent concerning  clones on infinite sets.)

  We write $\oo n $ for the set $ X^{X^n}$
of all $n$-ary operations. For a clone $C$, call  $C^{(n)}:= C
\cap \oo n $ the ``$n$-ary fragment of $C$''.

 $C^{(1)}$ is a submonoid of the monoid $X^X$ of all unary operations. For any monoid $M \subseteq X^X$ the set of all clones $C$ with $C^{(1)} = M$ is called the monoidal interval of $M$; it has a least element, the clone generated by $M$, and a largest element $\Pol(M)$, the set of all operations $f$ satisfying $f(m_1,\ldots, m_k)\in M$
 whenever $m_1,\ldots, m_k\in M$.  (Here $f(m_1,\ldots, m_k)$ is the unary operation
 mapping $x$ to $f(m_1(x), \ldots, m_k(x))$.

 In \cite{analytic}, we showed that on $X=\NN$ there are uncountably many clones containing
 all unary operations (but only two coatoms, see \cite{Gavrilov:1965}, \cite{737});  in other words, the monoidal
interval of $X^X$ is uncountable.
Pinsker in \cite{Pin06Monoidal} has constructed (on arbitrary infinite
base sets $X$) different monoids whose monoidal interval have various sizes,
among them also one  whose monoidal interval has size $2^{2^{|X|}}$).

We will show here that (for $|X|=\aleph_0$)
the interval associated with the monoid of
all $f:X\to X$ has
the largest possible size: $2^{\mathfrak c}$.   We will also
construct, for any natural number $n\ge 1$,  many clones which share their
$n$-ary fragment with $2^{\mathfrak c}$ other clones.

\subsection*{Local clones}

A clone $C$ is local if it is closed in the natural topology on $\mathscr O$,
that is:  for all $f\in \oo k  \setminus C$ there is a finite
set $A\subseteq X^k$ such that no $g\in C$ agrees
with $f$ on $A$.

For any infinite set $X$, the local clones on $X$ which contain all unary
operations are precisely known: they form an increasing chain of order type
$\omega + 1$.

\subsection{Main results}
\begin{Theorem}\label{theorem.1}
Let $X=\NN$ be countably infinite.  Then there are $2^{\frak c}$
clones on $X$ containing the monoid of all unary operations.
\end{Theorem}

To generalize the theorem also to larger arities, we need the
following technical definition:
\begin{Definition}  Let $\alpha\in \RR$.
  A operation $f:X^d \to X$ is defined to be $\alpha$-modest iff for all natural numbers
  $N$ and all $Y \subseteq X$ of cardinality $N$, the range of $f\on
  Y^d$ has at most $\alpha N$ elements.

  $f$ is modest if $f$ is $\alpha$-modest for some $\alpha$.

  We call a clone $C$  modest if all operations in $C$ are modest.

  We write $\mathscr M$ for the set of all modest operations.

\end{Definition}

Note that $\mathscr M$ is a clone (the greatest modest clone),
and that all unary operations are modest; in addition, all
operations with finite range are modest, as well.

\begin{Theorem}\label{theorem.2}
  Let $d\ge 1$, and let $C$ be a modest clone on $\NN$
  containing all $d$-ary
  operations with range $\{0,1\}$. Then there are $2^{\mathfrak c}$
  many clones $D$ with $D\cap \oo d = C \cap \oo d $.
\end{Theorem}
Taking $d=1$ and $C$ the clone of all essentially unary operations,
we get theorem~\ref{theorem.1} as a special case.

Machida \cite{Mac98} has defined a natural metric on clones:  The
distance between two clones is $1/n$, where $n$ is minimal with
$C\cap \oo n \not= D\cap \oo n $. In this language,
theorem~\ref{theorem.2} says that certain sets of clones can be
arbitrarily small from the metric/topological point of view, and
still large when measured by cardinality.

% \subsection{Notation}
Let $F$ be a set of operations. We write $\langle F \rangle $ for the smallest
clone containing $F$.   If $C$ is a clone, we may write $\langle F\rangle_C$
instead of $\langle F\cup C\rangle$ (and for $F=\{f,g,\ldots\}$ we
write $\langle f,g,\ldots\rangle_C$
instead of  $\langle\{ f,g,\ldots\}\rangle_C$.
Note that $f\in \langle F\rangle _C$ iff
there is a finite subset $F_0\subseteq F$ with
$f\in \langle F_0\rangle _C$.

Both sections of this paper use the following easy fact:
\begin{Lemma}\label{lemma.basic}
  Let $C$ be a clone, and let $(f_i:i\in I)$ be a family of operations
  which is independent over $C$
  (which means that $f_i \notin \langle f_j: j\not=i \rangle _C$ for all $i\in I$).
  For $J \subseteq I$ let $C_{J} = \langle f_i:i\in J\rangle_C$.
  Then:
  \begin{enumerate}
  \item[(a)] The map  $J \mapsto C_{J}$ is a 1-1
    order-preserving map from $\mathfrak P(I)$, the power set of $I$,
    into the interval~$[C, \langle
      f_i:i\in I\rangle]$ in the clone lattice (both ordered by inclusion).
  \item[(b)] If $I$ has cardinality $\kappa$, then $\{C_J:J \subseteq
    I\}$ contains   $2^\kappa $ many elements, and it is
    order-isomorphic with $\mathfrak P(I)$.
  \item[(c)] Assume moreover that
    $\{f_i:i\in I\} \subseteq \Pol(C\cap \oo d)$.
    Then $C_J \cap \oo d = C\cap \oo d $ for all $J \subseteq I$.
  \end{enumerate}
\end{Lemma}
\begin{proof}
  (a) and (b) are clear.

  For (c), recall that $\Pol(C\cap \oo d)$ is the set of all
  operations $f$ with $f(c_1,\ldots, c_m) \in C \cap \oo d $ whenever
  $c_1,\ldots, c_m\in C \cap \oo d $.  Clearly the assumption implies
  $$ C \subseteq \langle f_i:i\in I\rangle_{C} \subseteq \Pol(C \cap \oo d ),$$
  and by definition, the clones $C$ and $\Pol(C\cap \oo d)$ have the same
  $d$-ary fragment $D$. Consequently, the $d$-ary fragment of $C_{J}$
  is $D$, as well.
\end{proof}

\section{Sparse graphs and modest operations}

\begin{Definition}
  Let $(V,E)$ be a graph (i.e., $E \subseteq [V]^2$, where $[V]^2$ is the set of
   2-element subsets of $V$).

   We say that  $(V,E)$ is $(k,l)$-sparse, if for every
   $U\subseteq V$ of size at most $k$, the induced subgraph on $U$ has at
   most $l$ edges.
\end{Definition}

%%%%%%%%%%%%% EXPLANATORY PHARAGRAPH %%%%%%%%%%%%%%%%%%%

In order to help the reader, in this paragraph we are providing
a brief and informal explanation for the technical details of
the rest of this section. In Lemma \ref{lemma.basic.sparse} below we will show, that
for large enough $N$ and $0 < \varepsilon < \frac 12$, there
exist graphs ${\mathcal G}$ on $N$ vertices, whose $M$-sized subgraphs
(for small $M$) are sparse, while at the same time, these ${\mathcal G}$ have ``many''
edges: the number of their edges is at least
$N^{1+\varepsilon}$. Using this lemma, we will be able to
construct functions on finite domains, which have large range,
but the range of their restrictions to small sets remain small;
for the details see Lemma \ref{sparse.2.modest}. ``Gluing together" carefully an
infinite sequence of such operations, we obtain a set $S$ of
operations on $\NN$ such that $S$ is independent and has
cardinality ${\mathfrak c}$. Combining this with Lemma \ref{lemma.basic}, the proof of
Theorem \ref{theorem.1} will follow quickly.

%%%%%%%%%%%%% end of EXPLANATORY PARAGRAPH %%%%%%%%%%%%%

\begin{Definition}
  Let $M$, $N$ be natural numbers, $0<\e<\frac 12$.
  We write $M \ll_\e N$ if $M\cdot N^{2\e-1}< 1/10$.
\end{Definition}

\begin{Lemma}\label{lemma.basic.sparse}
  Let $\e < 1/2$, and let
  $1 \leq M \ll_\e N$.   Then there is a graph $G=(V,E)$ with $N$ vertices
  and more than $N^{1+\e} $ edges which is $(k,2k)$-sparse for all
  $k\le M$.
\end{Lemma}

\begin{proof}
   We will use an Erd\H os type probability argument: we will define a
   suitable probability measure on all graphs on $N$ vertices and then
   show that the set of graphs not satisfying the conclusion has small
   measure. \\
   \indent We note, that a somewhat stronger form of the lemma
   follows quickly from the Central Limit Theorem. For completeness, we present an
   elementary proof.

   \newcommand{\fn}{\frac{N(N-1)}{2}}

   \newcommand{\muone}{4N^{-1+\e}}
   Let $p:= \muone$,  and let
   $\mu$ be the probability measure on $\{0,1\}$ with
   $\mu(\{1\}) = p $.   Fix a set $V$ of $N$ vertices; there are
   $\fn $ potential edges.  Via characteristic functions, we
   identify the set of all graphs on $V$
   with the product space $ \{0,1\}^{\fn}$, which is equipped with the product probability structure.  In order
   to keep notation simpler, the product measure will also be called $\mu$.

   In other words, for each potential edge $e$ we flip a weighted coin
   (independent of all other coin flips),
   and with probability $p$ will decide to add $e$ to our graph.
   The expected number of edges is $\fn\cdot p  \approx
    2 N^{1+\e}$, with variance $\fn p(1-p) \approx 2
    N^{1+\e}$.  By Chebyshev's inequality, most graphs will have more
    than $N^{1+\e}$ edges.  More precisely, the measure of the set of
    graphs with fewer than $N^{1+\e}$  edges is smaller than
    \[ \frac{
      \fn p(1-p)
    }{
      (\fn\cdot p - N^{1+\e})^2
    }
    \approx
    \frac{
          2 N^{1+\e}
    }{
      (N^{1+\e})^2
    } =
 2
    N^{-1-\e}  < 1/2 , \]

   \noindent because, by the assumptions of the lemma, we have $4 \leq
   N$.\\
   \indent
   We now estimate the measure of the set ${\mathscr G}$ of all graphs on $V$
   which
   are not $(k,2k)$-sparse  for some $k\le M$.

   For any set $E'\subseteq [V]^2$ we let
   ${\mathscr G}_{E'}$ be the set of all graphs whose edges include the set
   $E'$.   Clearly $\mu ( {\mathscr G}_{E'}) = (\muone)^{|E'|}$.

   For each graph $(V,E)$ which is not $(k,2k)$-sparse there exists a set $V'$ of $k$
   vertices, and a set $E' \subseteq [V']^2$ with $2k$ elements such
   that  $E \supseteq E'$, i.e., $(V,E)\in {\mathscr G}_{E'}$.
      So the measure of all those graphs is bounded
   above by
   $$ % \sum_{k=5}^{k^2}
   \sum_{V'\subseteq V\atop |V'|=k}
   \sum_{E'\subseteq [V']^{2}\atop |E'|=2k} \mu ( {\mathscr
   G}_{E'}). $$ The crucial component in this sum is the summation over
   all subsets of size $k$; this will be estimated by a factor $N^k$;
   the other summations will be replaced by factors that depend on $k$
   only. Altogether we get an upper bound \\
   \\
   \centerline{$ N^{k}(k^{2})^{2k}(4N^{-1+\e})^{2k} =
   (2k)^{4k}N^{k}N^{-2k(1-\e)} = (2k)^{4k}N^{k(2\e-1)} \approx N^{k(2\e-1)} $.}  \\
   \\
   Now summing over all $k\le
   M$ yields $$ \sum_{k=1}^{M} N^{k(2\e-1)}
   \le M\cdot N^{2\e-1}  < 1/10$$ as $M\ll_\e N$.

   Hence the set of graphs satisfying the conclusion has measure  $>0$,
   so it is nonempty.
\end{proof}

\begin{Lemma} \label{sparse.123}
  Let $0<\e < \frac 12$.
  There is an increasing sequence $\langle N_\ell:\ell\in \NN\rangle$ of natural numbers and a
  sequence $\langle (V_\ell,E_\ell): \ell\in \NN \rangle$ of graphs such that:
    \begin{enumerate}
      \item $max \{N_{\ell}^{2}+1,2^{3N_{\ell}}, 1+|E_{\ell}|
      \} < N_{\ell +1}$.
      \item $V_\ell = [N_{\ell-1} , N_\ell)$.
      \item $|E_\ell|\ge N_\ell^{1+\e}$.
      \item For all $k\le 2^{\ell+1}N_{\ell-1}$, the graph $(V_\ell,E_\ell)$
is $(k,2k)$-sparse.
    \end{enumerate}
\end{Lemma}
\begin{proof}
   We can choose $N_\ell$ by recursion; given $N_{\ell-1}$,
   Lemma~\ref{lemma.basic.sparse} tells us how large $N_\ell$ has to
   be. In more detail, let $\e'$ be such that $\e < \e' < \frac 12$. Then,
   by Lemma
   \ref{lemma.basic.sparse}, there exist $N_{\ell}'$ and a graph
   ${\mathcal G}$ with $N_{\ell}'$ vertices and more than $(N_{\ell}')^{1+\e'}$
   edges which is $(k,2k)$-sparse for all $k \leq 2^{\ell + 1}N_{\ell -1}$.
   Enlarging $N_{\ell}'$ if necessary, we may assume, that 
\begin{itemize}
\item 
     $(1)$ holds (more precisely, $N_{\ell}'$ is larger
   than the left hand side of the $(\ell-1)^{th}$ instance of (1)), and 
   \item  $(1+\e)\ln(2) < (\e'-\e)\ln(N_{\ell}')$ and $2 N_{\ell - 1} \leq
   N_{\ell}'$.
\end{itemize}
   Take $N_{\ell} := N_{\ell-1} + N_{\ell}'$.
   Let ${\mathcal G}_{\ell}$ be an isomorphic copy of ${\mathcal G}$ with
   $V_\ell = [N_{\ell-1},N_{\ell})$. Now (2) and (4) of
   the statement clearly hold for ${\mathcal G}_{\ell}$. To check (3), it is enough to show that $N_{\ell}^{1+\e}
   \leq (N_{\ell}')^{1+\e'}$, that is, $$
   (*)\qquad \qquad \ln(N_{\ell}^{1+\e})
   \leq \ln((N_{\ell}')^{1+\e'}).$$
   The following calculation proves $(*)$: 
  \begin{align*}  
   \ln(N_{\ell}^{1+\e}) &=  (1+\e)\ln(N_{\ell - 1} +
   N_{\ell}') \\
   & \leq (1+\e)\ln(2N_{\ell}') \\
   &= (1+\e)\ln(N_{\ell}') +
   (1+\e)\ln(2) \\
   & \leq (1+\e)\ln(N_{\ell}') + (\e'-\e)\ln(N_{\ell}') \\
    &= (1+\e')\ln(N_{\ell}') \\
    &= \ln((N_{\ell}')^{1+\e})\qedhere
  \end{align*}
\end{proof}
So our graphs $(V_\ell,E_\ell)$ have ``many edges'' on a large
scale (i.e., looking at the whole graph), but only ``few edges''
on the small scale (looking at small induced subgraphs).

\begin{Definition}
  A $\varrho$-ary (partial) function $f:V^{\varrho}\to \NN$ is defined to be $(k,l)$-modest iff for any $U_0,\ldots, U_{\varrho-1} \subseteq V $ of size at most $ k$,
$f\on(
  U_0 \times \cdots  \times U_{\varrho-1})$ has at most $l$ values.

  The ``support'' of a function $f$ is the set of elements
  % $x$
  in the
  domain of $f$ where $f$ is defined and its value is not equal to~$0$.
\end{Definition}

\begin{Lemma}\label{sparse.2.modest}
  Let $(V,E)$ be a graph which is $(k,2k)$-sparse for all $k\le M$.
  Then there is a function $f:V\times V\to \NN$ which has at least $ |E|$
  values but is $(k,5k)$-modest for all $k\le M/2$.
\end{Lemma}
\begin{proof}  Let $f$ be a symmetric function which takes different
  values on all edges in $E$, and is constantly zero outside $E$. Then
  for each $U_1,U_2\subseteq V$ of size $k\le M/2$, $E\on (U_1\cup
  U_2)^2$ has at most $2\cdot 2k$ edges, so  $f$ can take at
  most  $4k+1 $ values on $U_1\times U_2 \subseteq (U_1\cup U_2)^2$.
\end{proof}

\begin{Corollary}\label{modest.123}
  There is an increasing sequence $\langle N_\ell:\ell\in \NN \rangle$ of natural numbers and a
  sequence $\langle s_\ell:\ell\in \NN \rangle$ of operations
  $s_\ell:[N_{\ell-1},N_\ell)^2 \to \NN$ satisfying the following:
    \begin{enumerate}
      \item $N_\ell^2 + 1 < N_{\ell+1}$ and also $2^{N_\ell} <
                  \sqrt[3]{ N_{\ell+1}}$.\label{modest.0}
      \item Each $s_\ell$ is $(k,5k)$-modest \label{modest.1}
    for all $k\le 2^\ell N_{\ell-1}$.
      \item Each $s_\ell$ is $(k,5k)$-modest for $k\ge N_{\ell+1}$.
    \label{modest.2}
      \item For all $\ell$, the range of $s_\ell$ has more than
    $N_\ell^{4/3}$ elements.
    \label{modest.3}
    \end{enumerate}
\end{Corollary}

\begin{proof} Let $\e = \frac 13$ and let $\langle N_{\ell}: \ell
\in \NN \rangle $ and $\langle (V_{\ell},E_{\ell}): \ell \in
\NN \rangle$ be the sequences obtained from
Lemma~\ref{sparse.123}. In addition, for every $\ell \in \NN$,
let $s_{\ell}$ be the operation obtained from
$(V_{\ell},E_{\ell})$ by Lemma~\ref{sparse.2.modest}. We claim
that this choice satisfies the statement. \\
\indent (1) follows from Lemma \ref{sparse.123}(1). Combining
Lemma \ref{sparse.123}(4) with Lemma \ref{sparse.2.modest} one
obtains (2). By the construction described in Lemma
\ref{sparse.2.modest}, the range of $s_{\ell}$ has cardinality at
most $|E_{\ell}|+1< N_{\ell+1}$. Hence (3) holds trivially because of Lemma
\ref{sparse.123}(1). Finally, (4) follows from Lemma
\ref{sparse.123}(3) (combined with the choice of $\e$ and with
Lemma \ref{sparse.2.modest}).

% { \tiny Property \ref{modest.0} together with Property
% \ref{fact:more}(\ref{small.range}) will trivially imply
% (\ref{modest.2}).
%
% The crucial point is the combination of (\ref{modest.1}) and
% (\ref{modest.3}).}

\end{proof}

\relax From now on we fix sequences $\langle N_\ell: \ell \in \NN \rangle$ and
$\langle s_\ell: \ell \in \NN \rangle$ as above.

\begin{Definition}
 For every $A\subseteq \NN $ let $s_A:\NN\times \NN\to
 \NN$ be defined from $s_\ell$ as follows: $s_A$ is $\bigcup_{\ell\in A}s_\ell$,
 extended by the value $0$ wherever it is undefined (i.e., $s_A\on
 [N_{\ell-1} , N_\ell) \times  [N_{i-1} , N_i)$ is constantly zero for
 $\ell\not=i $).
\end{Definition}

\begin{Lemma} \label{lemma:i.ell.A}\
  \begin{enumerate}
    \item  \label{i.not.ell}
      If $\ell < i $, then $s_i$ is
      $(k,5k)$-modest
      for all $k \le 2^\ell N_\ell$.
    \item \label{ell.not.A}
      If $\ell\notin A$, then $s_A$ is $(k,6k)$-modest
      for all $k$ in $[N_\ell, 2^\ell N_\ell]$.
  \end{enumerate}
\end{Lemma}
\begin{proof}
  First we prove (\ref{i.not.ell}).
  By Lemma
  \ref{modest.123}(\ref{modest.1}), $s_i$ is
  $(k,5k)$-modest for all $k\le 2^i N_{i-1}$, so certainly also
  for all $k\le 2^\ell N_\ell$.

  Now  we prove
  (\ref{ell.not.A}).   Let $X,Y $ be sets of size $k$,
  with  $k$ in $[N_\ell, 2^\ell N_\ell]$.  Let $X_-=X\cap N_\ell$,
  $X_+=X\setminus X_-$, and define $Y_-,Y_+$ similarly.
  We have $$s_A[ X \times Y ] \subseteq \ \
  s_A[X_-\times Y_-] \ \cup  \  s_A[X_+\times Y_+]\  \cup \ \{0\}.$$
  \begin{itemize}
  \item  Because $\ell \not \in A$, $s_A$ is constantly 0
    on $(X_-\times Y_-) \setminus (N_{\ell-1} \times N_{\ell-1})$. Hence the
    first set has size at most $N_{\ell-1}^2\le N_\ell-1 \le k-1$.
  \item To estimate the size of $s_A[X_+\times Y_+] \cup \{0\}
  \subseteq  \{0\} \cup \bigcup_{i>\ell} s_i[X_+\times Y_+]$ we use
  (\ref{i.not.ell}) to see that this set is bounded by $5k+1$.
  \item So $s_A[ X \times Y ]$ has size at most $6k$.
  \end{itemize}

\end{proof}

\begin{Definition}
  Let $A_1,\ldots, A_n \subseteq \NN$.  A (binary) \emph{term} in the
  operations $s_{A_1},\ldots, s_{A_n}$ is a formal expression involving
  (some of) the variables $x,y$, (some of) the
  operations $s_{A_1},\ldots, s_{A_n}$, as well as any unary operations.
  (We trust the reader to supply a formal definition by induction.)

  The \emph{depth} of a term $\tau$ is defined inductively as follows:
  \begin{itemize}
    \item $x$ and $y$ have depth $0$.
    \item For any unary operation $u$, the depth of $u(\tau)$ is $1$ more
     that  the depth of $\tau$.
    \item Let $m$ be the maximum of the depths of $\tau_1$ and
      $\tau_2$. Then the depth of $f_{A_i}(\tau_1,\tau_2)$ is  $m+1$.
  \end{itemize}
\end{Definition}

Every term naturally induces a binary operation on $\NN$. (Note
that the same operation may be represented by different terms,
even terms of different depths.)
% For example, the operation represented by the term
% $f_{A_1}(y,f_{A_2}(y))$ (of depth~2)
%
% Here f_{A_{2}} should be a binary operation !!!
%
% can also be represented by a term of depth~0.

\begin{Lemma}\label{lemma:main} Let $\tau$ be a term in the
  operations $s_{A_1},\ldots, s_{A_n}$  of depth $d$.    Let $\ell>d\log_{2}(6)$,
  and assume $\ell\notin A_1 \cup \cdots \cup A_n$.  Then we have:
  \begin{enumerate}
  \item [(1)]  The operation
  represented by $\tau$ is $(N_\ell,6^d N_{\ell} )$-modest.
  \item [(2)] In particular, $\tau $ cannot represent the operation $s_\ell$, or
  $s_B$ for any $B$ containing $\ell$.
  \end{enumerate}
\end{Lemma}

\begin{proof}
% [Proof of Lemma~\ref{lemma:main}]
  We start to show (1) by induction on $d $ (or more precisely, on $\tau$).
  \begin{itemize}
    \item If $\tau $ is $x$ or $y$, this is trivial.
    \item If $\tau = u(\tau_1)$ then again the range of $u(\tau_1)$ is
      not larger than the range of $\tau_1$.
   \item Assume $\tau = s_{A_i}(\tau_1,\tau_2)$, where the depths of $\tau_1$ and
   $\tau_2$ are at most~$d$. Observe that
     \begin{itemize}
       \item Both $\tau_1$ and $\tau_2$ are
     $(N_\ell,6^d N_\ell)$-modest by the inductive assumption.
       \item
     By Lemma \ref{lemma:i.ell.A}(2), $f_{A_i}$ is
     $(6^d N_\ell, 6\cdot 6^{d} N_\ell)$-modest.
     (Recall that $d\log_{2}(6)\le \ell$,
     so $6^d N_\ell \le 2^\ell N_\ell$.)
       % \item So by
     % \ref{fact:more}(\ref{composition}) we get that
     % $\tau$ is $(N_\ell,6^{d+1} N_\ell)$-modest.
     \end{itemize}
     Now let $U_{1},U_{2} \subseteq \NN$ be two sets, both
     of size  at most
     $N_{\ell}$. Then, according to the previous
     observation, the ranges of $\tau_{1}\on U_{1}$ and
     $\tau_{2}\on U_{2}$ have size at most $6^{d}N_{\ell}$. Hence,
     again by the previous observation, the cardinality of the range of $\tau \on U_{1} \times
     U_{2}$ is at most $6 \cdot 6^{d} N_{\ell} = 6^{d+1}
     N_{\ell}$, as desired. \\
     \indent Now we turn to prove (2). By assumption, $ 6^d \le 2^{\ell} \le 2^{N_{\ell-1}} $,
  so by (1) of Corollary \ref{modest.123} we have $6^d N_\ell <
  N_\ell^{4/3}$. Hence, by (1) of the present theorem,
  $|ran(\tau\on N_{\ell} \times N_{\ell}) | \leq 6^d N_\ell < N_{\ell}^{4/3}$, while,
  according to Corollary \ref{modest.123} (4), we have
  $|ran(s_{B} \on N_{\ell})| > N_{\ell}^{4/3}$.
  \end{itemize}
\end{proof}

\begin{Corollary} \label{coroll:indep}
  Let $B, A_1,\ldots, A_n$ be pairwise distinct subsets of $\NN$ such that
  $B \setminus  (A_1 \cup \cdots \cup A_n)$ is infinite.
  Then
  $s_B\notin \langle s_{A_1},\ldots, s_{A_n}\rangle_{\oo 1 }$.
\end{Corollary}

\begin{proof}
Assume, seeking a contradiction, that $s_B\in \langle
s_{A_1},\ldots, s_{A_n}\rangle_{\oo 1 }$. Then there exists a
term $\tau$ in $A_{1}, \ldots A_{n}$ representing $s_{B}$. Let
$d$ be the depth of $\tau$. Then there exists $\ell \in B
\setminus  (A_1 \cup \cdots \cup A_n)$ with $\ell > d\log_{2}(6)$.
Then, by Lemma \ref{lemma:main}(2), $\tau$ does not represent
$s_{B}$. This contradiction completes the proof.
\end{proof}

\begin{Fact} There exists an independent family $(A_r:r\in \mathbb R)$
of $\mathfrak c$ subsets of~$\NN$.
That is, for all disjoint finite
 subsets $J_+,J_-\subseteq \mathbb R$ the set $\bigcap_{r\in I_+} A_i \ \cap \
 \bigcap_{r\in I_-} (\NN\setminus A_i) $  is nonempty and even
 infinite.
\end{Fact}
\begin{proof} This is well known.  For example, replacing the base set
  $\NN$ by $\mathbb Q[x]$,  the set of all polynomials with
  rational coefficients, we can take $A_r:= \{p(x)\in \mathbb Q[x]:
  p(r)>0\}$.
\end{proof}

\begin{proof}[Proof of theorem \ref{theorem.1}]
 % Let $R$ be a set of size~$2^{\aleph_0}$.   We can find an
 Choose an independent family $(A_r:r\in \mathbb R)$ of subsets of $\NN $.
 Then, for all finite $S \subseteq \mathbb R$ and all
 $r\in\mathbb  R\setminus S$ the set $A_r\setminus \bigcup_{s\in S} A_s$ is
 infinite.
 % Let $C$ be the clone generated by the set of all unary operations on $\NN$.
 By Corollary \ref{coroll:indep}, $\{ s_{A_{r}} :r\in \mathbb R
 \}$ is a family of operations independent over $\oo 1$: for any $r \in
 \mathbb R$, we have $s_{A_r} \not \in \langle s_{A_{p}} :p
\in \mathbb R \setminus \{ r
 \} \rangle_{M \cap \oo d}$.
By Lemma \ref{lemma.basic} we are done.

\end{proof}

\section{Higher arities}
\begin{Definition}
We say that an operation $f:X^d\to X$ is \emph{modest} if there is some
$k$ such that for all $N>1$, $f$ is $(N,kN)$-modest, i.e.:
\begin{quote}
$f[X_1\times \cdots \times X_d] $ has at most $kN$ elements
\\
whenever each set  $X_i \subseteq X $ has at most $N$ elements.
\end{quote}

We call a clone $C$  modest if all operations in $C$ are modest.
\end{Definition}

Note that the set of all modest operations is a clone (the
greatest modest clone), and that all unary operations are modest,
as are all operations with finite range. \\
\indent The main result of the present section is as follows.

\begin{Theorem}\label{theorem.2.again}
  Let $d\ge 1$, and let $C$ be a modest clone containing all $d$-ary
  operations with range~$\{0,1\}$. Then there are $2^{\mathfrak c}$
  many clones $D$ with $D\cap \oo d = C \cap \oo d $.
\end{Theorem}

We postpone the proof of this theorem to the end of this section.
The number $d$ will be fixed through this section.

\begin{Definition}  \label{def:term1} \
\begin{itemize}
\item
 We fix a language with object variables
  $x_i$, $i\in \NN $ and formal operation variables ${\tt f}^i_j$
  ($i,j\in \NN $), where the superscript $i$ denotes the formal
  arity of ${\tt f}^i_j$. Terms are defined as usual: each object
  variable is a term, and whenever $t_1,\ldots, t_i$ are terms and
  $j\in \NN$, then ${\tt f}^i_j(t_1,\ldots, t_i)$ is a term, as well.
\item   The set of all terms can be enumerated as
  $\{\tau_1,\tau_2,\ldots\}$ such that
  % the term
  $\tau_m$ contains at
  most $m$ occurrences of operation symbols, and each operation
  symbol occurring in $\tau_{m}$ is at most $m$-ary.
\item   Let $\tau $ be  a term.  We say that a sequence of
functions $\bar g =
  (g^i_j:(i,j)\in S)$ is suitable for $\tau$ iff each $g^i_j$ has arity $i$, and
  $(i,j)\in S$ whenever the variable ${\tt f}^i_j$ appears in~$\tau$.
\item   If $\bar g$ is suitable for $\tau$, then
  plugging in the $g^i_j$ for the ${\tt f}^i_j$ will
  yield a $d$-ary operation on $X$ which we denote by $\tau[\bar g]$.
\end{itemize}
\end{Definition}

\begin{Definition}
  Let $d\ge 2$. For any set $V$ we let $[V]^d$ be the set of
  $d$-element subsets of $V$.
  The structure $(V,E)$ is defined to be a
  %  (or sometimes just the set $E$, if $V$ can be
  %   deduced from the context)
  \emph{$d$-uniform hypergraph} iff $E \subseteq [V]^d$.  The elements of~$E$ are called
  the \emph{ hyperedges} of~$(V,E)$.

   Every $V' \subseteq V$ naturally induces a hypergraph $(V', E \cap
   [V']^d)$, which we may also denote by $(V', E\on V')$.

   We say that  $(V,E)$ is $(k,l)$-sparse, iff for every
   $X\subseteq V$ of size at most $k$, the hypergraph $(X,E\on X) $ has at
   most $l$ hyperedges.
\end{Definition}

\begin{Lemma}\label{FUNCTIONS}
  Fix $d$, $k$, $\e$.
  Let $V$ be a set of cardinality $N$,
  % respectively, with $V  \subseteq  W$. % and $k\le \log N$
  let $(V,E)$ be a
  $(d+1)$-uniform  hypergraph with at least $N^{d+\e}$ hyperedges.
  % Moreover, assume $5 \log N < N^\e$.
  %
  % Then
  If $N$ is large enough, then
  there is an operation $s:V^{d+1} \to V$ whose support is contained in $E$ and whose values
  are in $\{0,1\}$ such that for any set $W$ with $V \subseteq W, |W| \leq kN$ the following holds.
  \begin{quote}
    Whenever $\tau\in \{\tau_1,\ldots, \tau_k\}$, \\
    $\bar{ \bar{g}}  = (g^i_j)_{i,j}$ is a suitable sequence of operations
    for $\tau$
    on $W$ with each $g^i_j$ being
    \begin{itemize}
    \item either of arity at most $ d$
    \item or of arity
      $d+1$ with support of size at most $ 3N\log N$
    \end{itemize}
    then $\tau[\bar{g}]$ does not represent~$s$. In particular, there exists $e \in E$ such that $s$ and $\tau[\bar{g}]$ have different values on $e$.
\end{quote}

\smallskip
\noindent
If $N$ satisfies the above conditions, then we will say that  $N$ is $k$-large.

\end{Lemma}
\begin{proof}  Let $W$ be a set containing $V$ with $|W| = kN$. Clearly, it is
enough to show that there exists an operation $s:V^{d+1} \rightarrow V$
satisfying the statement for this particular $W$. There are only
$(kN)^{(kN)^{d}}$ $d$-ary operations on~$W$, and
  only $k$ terms to be considered.  A support is a subset of
  $[W]^{d+1}$;
  % since we only have to consider supports of size $2N$,
  there
  are fewer than
  $\binom{(kN)^{d+1}}{3N\log N} \le (kN)^{3N\log (N)(d+1)}$ possible supports of size
~$3N\log N$. For each support of size $3N\log N$ there are at most
$(kN)^{3N\log N}$
  possible operations that have this support. By the enumeration
  fixed in Definition \ref{def:term1}, each term $\tau_{i} \ (i \leq
  k)$ contains at most $k$ many operation variables. Counting the possibilities of choosing $k$
  many $d$-ary operations and $k$ many $d+1$-ary operations with
  support of size at most $3N\log N$, one can see, that altogether there are fewer than \\
  \\
  \centerline{$t:= (kN)^{(kN)^d \cdot k} \cdot
       k\cdot (kN)^{3N\log(N)(d+1)k} \cdot (kN)^{3N\log(N)k}$} \\
       \\
  % (N^2)^{2N}\le N^{ 1+2d+2+4N}
  % \]
  operations represented by such terms. We may assume $k\le \log N$. Estimating $k$ by $N$ or
  by
  $\log N$, one obtains \\
  \\
  \centerline{$ t \leq (\log N) \cdot (N \cdot \log N)^{\log N \cdot (N \log N)^{d}} \cdot N^{6N(d+1)\log^{2} N} \cdot
        N^{6N \log^{2} N}.$} \\
  \\
  Recall, that for any $\delta > 0, \ d \in \NN$ and for large enough $N$, one has $\log^{d} N \leq N^{\delta}$.
  Let $0 < \delta  < \e $. Then for large enough $N$, each of the four factors of $t$ can be estimated by
  $N^{\frac 14 \cdot N^{d+\delta}}$. Consequently, for large enough $N$, we have $t < N^{ N^{d+\delta}} = 2^{N^{d+\delta}\cdot \log N}$.
  This number (for large enough $N$), is certainly less
  than
  $
  % N^{2N^d(\log^d N)} <2^{5 N^d\,\log N}$.
  % $N^{N^{d +\frac \e2}} <
  2^{N^{d + \e}}$.

  But there are at least $2^{N^{d+\e}}$ possible operations on $E$ with values
  in~$\{0,1\}$.  So not all of them are representable.
\end{proof}

\begin{Lemma}\label{GRAPHS}
  Let $0<\e < 1/2$.
  Then there are sequences $\bar N = \langle N_\ell:\ell< \NN \rangle$,
  $\bar E=\langle E_\ell:\ell< \NN \rangle $
  with the following properties:
  \begin{enumerate}
  \item $\bar N$ is strictly increasing and in fact $N_{\ell-1}^{d+1} <
    N_\ell$, $2^{\ell} \leq N_{\ell}$
    and $N_{\ell}$ is ${\ell}$-large for all~$\ell$.   We will write $V_\ell$ for the interval
    $[N_{\ell-1} , N_\ell)$.
  \item $(V_\ell, E_\ell)$ is a $(d+1)$-uniform hypergraph with more
    than $N_\ell^{d+\e}$ hyperedges.
  \item For every $k \le N_{\ell-1} ^2$, $(V_\ell, E_\ell)$
    is $(k,2k)$-sparse.
  \end{enumerate}
\end{Lemma}
\begin{proof}  This proof is only a slight variation of the proof of
    Lemma~\ref{modest.123}, so we will be brief.

    Assume $N_{\ell-1}$ has already been defined. We will choose $N_{\ell}$ after a certain amount of extra work such that $N_\ell \gg N_{\ell-1}$.
    Assume for a moment, that $N_{\ell}$ is already defined. Let
    $V_\ell:= [N_{\ell-1}, N_\ell)$.  Let $J$ be the cardinality of
    the set $[V_\ell]^{d+1}$ of all potential hyperedges:
        $J = \binom{N_\ell-N_{\ell-1} }{d+1}$.

    On the set of all
    $(d+1)$-uniform hypergraphs (which we may identify with
    $2^J$), we define a product measure by declaring the probability
    of each potential hyperedge to be $p:= 2(d+1)!\cdot N^{\e-1}$.

    So the expected number of hyperedges of a random hypergraph
    is $p J =    2(d+1)!\cdot N_\ell^{\e-1} \cdot \binom{N_\ell-N_{\ell-1}
    }{d+1} \approx 2 N_\ell^{\e-1}\cdot N_\ell^{d+1} = 2 N_{\ell}^{d+\e}
    $. Again using Chebyshev's inequality, we see that with high probability a random
    hypergraph will have more than $N_{\ell}^{d+\e}$ hyperedges.

   Now we estimate the probability that there is a sub-hypergraph with
   $k\le N_{\ell-1}^2$ vertices which has more than $2k$ hyperedges,
   and we will show that it is very low.

   For each potential $k$ there are at most $ \binom{N_\ell}{k} \le N_\ell^k$
   subsets;
   for each such subset $S$, the probability that a given set $H$
   of hyperedges
   with $j :=|H|\ge 2k$ appears as a subset of $E\on S$
   is $\le p^j\le p^{2k}$.  There are
   $\binom{k^d}{j}\le 2^{k^d}$ possibilities for~$H$.
   So the probability that
   such a bad subgraph of size $k$ exists is bounded from above by
  \[
           N_\ell^k \cdot p^{2k} \cdot 2^{k^d}.
  \]
  There are $N_{\ell-1}^2$ possibilities for~$k$, so we have to
  choose $N_{\ell}$ such that $$ (*) \indent
  \sum_{k=1}^{N_{\ell-1}^{2}} N_{\ell}^{k} \cdot p^{2k}2^{k^{d}} \leq
  \frac 12. $$

   But $N_\ell^k \cdot p^{2k} \approx N_{\ell}^k  N_{\ell}^{(\e-1)2k} =
   N_{\ell}^{k(2\e-1)}$ which converges to $0$ if $N_{\ell}$
   converges to infinity. Hence, one may choose $N_{\ell}$ so
   large, that \\
   \\
   \centerline{$\displaystyle{N_\ell^k \cdot p^{2k} < \frac{1}{N_{\ell-1}^{2} \cdot
   2^{(N_{\ell-1}^{2})^{d}}}}$} \\
   \\
   and $N_{\ell} >  max \{2^{\ell},N_{\ell-1}^{d+1}\}$ hold. Further increasing $N_{\ell}$ if necessary, we may choose it to be $\ell$-large, as well. Estimating $2^{k^{d}}$ by
   $2^{(N_{\ell-1}^{2})^{d}}$ in the left hand side of $(*)$, it follows,
   that the inequality in $(*)$ holds. \\
   \indent So the set of hypergraphs on $V_{\ell}$ which are not $(k,2k)$-sparse for some $k \leq N_{\ell-1}^{2}$  has measure at most
   $\frac 12$, while, almost all hypergraphs on $V_{\ell}$ have
   $N_{\ell}^{d+\e}$ hyperedges. It follows, that there exist
   $N_{\ell}$ and $E_{\ell}$ satisfying the requirements of the lemma,
   and thus, the sequences in the statement can be constructed
   recursively.
\end{proof}

\begin{Definition}
\label{def:wl}
  Let $\bar N$ and $\bar E$ be as in Lemma \ref{GRAPHS}.
  % and let $\ell \geq 2$.
  For each $V_\ell= [N_{\ell-1},N_\ell)$
  % let $W_\ell $ be a set of size $\ell |V_\ell|$
  % containing $V_\ell \cup \{0,1\}$, and
  let $s_\ell$ be a $(d+1)$-ary operation
    with support $E_\ell$
    which differs on $E_\ell$ from
    each $\tau_i[g]$ ($i\le \ell$, $\bar g$ as in Lemma \ref{FUNCTIONS}).

  For each infinite $A \subseteq \NN$ let
  $s_A:=\bigcup_{\ell\in A} s_\ell$ (where we replace all undefined values of
  $s_A$ with $0$).
\end{Definition}

\begin{Lemma}
\label{lemma:support} Let $B \subseteq \NN$ be infinite, and assume
$\ell \in \NN \setminus B$. Let $W \subseteq \NN$ be such that
$|W| \leq \ell \cdot N_{\ell}$. Then the cardinality of the
support of $s_{B} \on W^{d+1}$ is at most $N_{\ell}(1+2\log
N_{\ell})$.
\end{Lemma}

\begin{proof} Throughout this proof, the support of a function $f$
is denoted by $supp(f)$.  Let $W_{1}=W \cap [0,N_{\ell-1}), \
W_{2} = W \cap [N_{\ell-1},N_{\ell}), \ W_{3} = W \setminus(W_{1}
\cup
W_{2})$. By construction, 
$$supp(s_{B} \on W^{d+1}) \subseteq supp(s_{B} \on
W_{1}^{d+1}) \cup supp(s_{B} \on W_{2}^{d+1}) \cup supp(s_{B} \on
W_{3}^{d+1})$$
\begin{itemize}
\item  Clearly, $|supp(s_{B} \on W_{1}^{d+1})| \leq
N_{\ell-1}^{d+1}$ and $N_{\ell-1}^{d+1} \leq N_{\ell}$ by Lemma
\ref{GRAPHS}(1). 
\item  In addition, $supp(s_{B} \on W_{2}^{d+1})$ is empty
because $\ell \not \in B$. 
\item  Clearly, $|W_{3}|\leq |W| \leq \ell \cdot N_{\ell} \leq
\log(N_{\ell})N_{\ell}$ (in the last estimation we used Lemma
\ref{GRAPHS} (1): $\ell \leq \log N_{\ell}$). In addition, by
Lemma \ref{GRAPHS} (3), for any $j > \ell$, $(V_{j},E_{j})$ is
$(N_{\ell}\log N_{\ell}, 2N_{\ell}\log N_{\ell})$-sparse. It
follows, that $ |supp(s_{B} \on W_{3}^{d+1})| \leq 2N_{\ell}\log
N_{\ell}$. 
\end{itemize}
Combining these observations, the statement follows.
\end{proof}

\begin{Lemma} \label{easy.alpha.beta}
 If $f_1,\ldots, f_m$ are $(k, k')$-modest $d$-ary
  operations,  $g$ is a $(k', k'')$-modest $m$-ary
  operation, then $g(f_1,\ldots, f_m)$ is $(k, k'')$-modest.
\end{Lemma}
\begin{proof}  Easy.
\end{proof}

\begin{Lemma}
\label{lemma:indep3}
 Let $\mathscr M$ be the clone of all modest operations.
 Let
 $A \setminus (B_1\cup \cdots B_r)$ be infinite.
 Then $s_A\notin \langle  s_{B_1}, \ldots, s_{B_r}\rangle_{M\cap \oo d}$.
\end{Lemma}

\begin{proof}
  For any term $\tau$ and any suitable sequence $\bar g$ (consisting
  only of operations in $ \langle (\mathscr M\cap \oo d ) \cup \{ s_{B_1},
  \ldots, s_{B_r }\}\rangle)$ we will find $\ell\in A$ such that
  $\tau[g]$ disagrees with $s_\ell$ (hence also with $s_A$)
  on~$E_\ell$.

  So fix a term $\tau = \tau_i$ and $\bar g$.  Let $\nu$ be the number of subterms of $\tau$ and let $k$ witness that all operations in $\bar g$
  are modest. Let $\ell >  \nu \cdot k^i $ be in
  $A \setminus (B_1\cup \cdots B_r)$. We claim that, for each subterm $\sigma$ of $\tau$ (of depth $s$), the range
  of $\sigma[\bar g]$ over the domain $V_\ell^{d+1}$  has cardinality
  at most $N_\ell \cdot k^s$.

  This can be proved by induction on the depth of $\sigma$, using Lemma \ref{easy.alpha.beta} combined with the fact that
  the operations $s_{B_j}$ take only 2 values, and that all other operations in
  $\bar{g}$ are modest, witnessed by $k$.

  Recall, that according to the enumeration fixed in Definition \ref{def:term1}, the depth of $\tau=\tau_{i}$ is at most $i$. So the set of all intermediate values in the computation of
  $\tau[g]$ on $E_\ell$ has size at most $ \nu \cdot k^i N_\ell < \ell N_\ell$. Let $W
  \supseteq V_\ell$ be a set of size at most $ \ell N_\ell$ containing $\{0,1\}$ and all
  these intermediate values. The term $\tau$ induces a partial
  function $\tau[\bar{g}]\on E_\ell$. By replacing all values of the
  operations in $\bar{g}$ by 0 if they are outside $W$, we get a sequence
  $\bar g'$ of operations with the following properties:
   \begin{itemize}
       \item $\tau[\bar g ']$ is a total function from $W^{d+1}$ to~$W$.
       \item $\tau[\bar g ']$ agrees with $\tau[\bar g]$ on~$E_\ell$.
       \item All operations   in $\bar g'$ are either some $s_{B_j}$
            or an operation of arity at most $ d$.
   \end{itemize}

   % Let $\bar g''$ be the sequence of functions obtained from $\bar g'$ by replacing each $s_{B_j}$ with the function mapping each element of $W$ onto $0$.
   %
   % Since $\ell \not\in B_1 \cup ... \cup B_r$, it follows, that for all $j<r$, the support of $s_{B_j} \on E_{\ell}$ is empty. Hence $\tau[\bar g''] \on E_{\ell} = \tau[\bar g' ] \on E_{\ell} = \tau[\bar{g}] \on E_{\ell}$.
   % The support of each $d+1$-ary function in $\bar{g''}$ is empty,
   %
   By Lemma \ref{lemma:support}, the support of each $s_{B_{j}}
   \on W^{d+1}$ is at most $N_{\ell}(1+2\log N_{\ell}) \leq 3
   N_{\ell}\log N_{\ell}$. So
   by the construction of $s_\ell$, and by Lemma \ref{FUNCTIONS}, $s_\ell$ disagrees with
   $\tau[\bar g']$ somewhere on $E_{\ell}$; so $s_\ell$ also disagrees
   with $\tau [\bar g]$.
\end{proof}

Now we are ready to prove Theorem \ref{theorem.2.again}.

%
%   Let $d\ge 1$, and let $C$ be a modest clone containing all $d$-ary
%   operations with range~$\{0,1\}$. Then there are $2^{\mathfrak c}$
%   many clones $D$ with $D\cap \oo d = C \cap \oo d $.
%
%

\begin{proof}[Proof of theorem \ref{theorem.2.again}]
 % Let $R$ be a set of size~$2^{\aleph_0}$.   We can find an
 Similarly to the proof of Theorem \ref{theorem.1},
 choose an independent family $(A_r:r\in \mathbb R)$ of subsets of $\NN $.
 Then, for all finite $S \subseteq \mathbb R$ and all
 $r\in\mathbb  R\setminus S$ the set $A_r\setminus \bigcup_{s\in S} A_s$ is
 infinite.
 % Let $C$ be the clone generated by the set of all unary operations on $\NN$.
 By Lemma \ref{lemma:indep3}, $\{ s_{A_{r}} :r\in \mathbb R
 \}$ is a family of operations independent over $M \cap \oo d$: for any $r \in
 \mathbb R$, we have $s_{A_r} \not \in \langle s_{A_{p}} :p\in \mathbb R \setminus \{ r
 \} \rangle_{\oo 1}$.
By Lemma \ref{lemma.basic} we are done.

\end{proof}

\begin{Corollary}
There exists a clone $C$ on $\NN$ such that for any $d \in \NN$ there are $2^{\mathfrak c}$ clones $D$ with
$C \cap \oo d = D \cap \oo d$.
\end{Corollary}

\begin{proof}
Let $C$ be the clone generated by all operations whose ranges are
a subset of $\{0,1\}$. To check, that this $C$ satisfies the
statement of the corollary, let $d \in \NN$ and let $C'$ be the
clone generated by all at most $d$-ary operations whose ranges
are contained in $\{0,1\}$. Then $C \cap \oo d = C' \cap \oo d$
and $C'$ is modest. Therefore, by Theorem \ref{theorem.2.again}
there exist $2^{\mathfrak c}$ many clones $D$ with $D \cap \oo d
= C' \cap \oo d = C \cap \oo d$.
\end{proof}

\bibliography{other,survey} \bibliographystyle{plain}

\end{document}